\documentclass{birkjour}
 \usepackage{graphicx}
 \usepackage{amsthm}
\usepackage{hyperref}
\usepackage{xcolor}
\newenvironment{refproof}[1]
  {\par\noindent\textbf{Proof of #1.}\quad}
  {\hfill$\square$\par}
\usepackage{algorithm}
\usepackage{algorithmic}
\usepackage{amsmath}
\usepackage{amssymb}
\usepackage{multirow}
\newcommand{\nc}{\mathbb{N}}

\usepackage{hhline}
\newtheorem{theorem}{Theorem}[section]
\newtheorem{defin}[theorem]{Definition}
\newtheorem{lemma}[theorem]{Lemma}

\newtheorem{cor}{Corollary}[section]
\newtheorem{rem}{Remark}[section]
\newtheorem{pro}{Proposition}[section]

\begin{document}

%-------------------------------------------------------------------------
% editorial commands: to be inserted by the editorial office
%
%\firstpage{1} \volume{228} \Copyrightyear{2004} \DOI{003-0001}
%\seriesextraline{This is the Concrete Title of this Book\br H.E. R and S.T.C. W, Eds.}
%\firstpage{1}

\title[On the Tameness of Power Series Space Pairs]{On the Tameness of Power Series Space Pairs}
\author{Buket Can Bahad\i r}

%----------Author 1

\address{Ankara, Turkey}
\email{ canbuket@gmail.com}

%\thanks{This work was completed with the support of our
%----------Author 2
%\author{A Second Author}
%\address{The address of\br
%the second author\br
%sitting somewhere\br
%in the world}
%\email{dont@know.who.knows}


%\date{January 1, 2004}

\begin{abstract}
				In this paper, it is shown that the tameness of the Köthe space pair $(\lambda^p(A),\lambda^q(B))$ is determined solely by the tameness of the family of quasi-diagonal operators defined between the pair of spaces. We use this tool to fill the gaps in characterization of pairs of power series spaces, adding to the previously established results of Dubinsky, Vogt \cite{vogt1}, Nyberg \cite{nyberg} and etc., and summarize this complete characterization in Table 1. As a result, we also show that the range of every continuous tame operator defined between power series spaces of infinite type has a basis. 
			\end{abstract}

%%% ----------------------------------------------------------------------
\maketitle
%%% ----------------------------------------------------------------------
%%%%%%%%%%%%%%%%%%%%%%%%%%%%%%%%%%%%%%%%%%%%%%%%%%%%%%%%%%%%%%%%%%%%%%%%%

\section{Introduction}
			The concept of tameness provides a tool to control the behavior of operators defined between Fr\'echet spaces. Tameness indicates that operators function within a well-defined structure and exhibit regular, non-pathological behavior. For this reason, it provides predictability and order while working with these operators in Fréchet spaces. 
			
			Tame spaces were first defined by Dubinsky and Vogt \cite{vogt1} in 1989, although the name appears in a paper of Hamilton written in 1982 \cite{hamil} regarding to linearly tame spaces. In his paper, Hamilton showed that it is possible to give an inverse function theorem for linearly tame spaces. Dubinsky and Vogt in \cite{vogt1} characterized the tame infinite type power series spaces. Vogt \cite{vogt3} showed that every finite type of power series space is tame. Nyberg \cite{nyberg} classified tame pairs of power series spaces. Later, Piszczek characterized tame pairs $(X,Y)$ in Fr\'echet spaces, where one of the spaces is a power series space, by using topological invariants.
			
			This concept can also be applied to the longstanding problem of determining whether every complemented subspace of Köthe spaces possesses a basis. Mitiagin and Henkin \cite{MH} provided an positive answer for all power series spaces of finite type. However, addressing this question for power series spaces of infinite type took considerably longer. In \cite{vogt1}, Dubinsky and Vogt showed that if a power series space of infinite type is tame, then the answer to this question is also positive. The most general solution to this problem was provided by Aytuna et al. in \cite{AKT}. They proved that a nuclear Fr\'echet space $E$ with properties $\underline{DN}$ and $\Omega$  contains a complemented copy of $\Lambda_{\infty}(\varepsilon)$,  provided that the diametral dimensions of E and $\Lambda_{\infty}(\varepsilon)$ are equal and $\varepsilon$ is stable. In \cite{A}, Aytuna  demonstrated that for a nuclear Fr\'echet space $E$ with properties $\underline{DN}$ and $\Omega$, and a stable finitely nuclear associated exponent sequence $\varepsilon$,
			$E$ is isomorphic to a power series space of finite type if and only if  $E$ is tame and 
			the approximate diametral dimensions of E and $\Lambda_{0}(\varepsilon)$ are equal. 
			Aytuna further utilized this result to characterize the tameness of spaces of analytic functions defined on Stein manifolds. For a detailed discussion on diametral dimension and approximate diametral dimension in nuclear spaces with the properties $\underline{DN}$ and $\Omega$, we refer to \cite{ND1}, \cite{ND2}.
			
			After establishing terminology and notation, in section 3, we showed that the tameness of Köthe space pairs
			$(\lambda^p(A),\lambda^q(B))$ can be approached by examining the tameness of only quasi-diagonal operators. In section 4, we looked at the results found in the literature, characterized the tameness of power series space pairs completely by filling the gaps in Table 1. In section 5, we proved that the range of every continuous tame operator defined between power series spaces of infinite type has a basis.
			
			\section{Preliminaries}
			
			In this section, we give some necessary definitions, terminology and notations; and present some fundamental results about tameness.
			
			A Fr\'echet space $E$ is a complete metric linear space whose topology is defined by an increasing sequence $(\|\cdot\|_{n})_{n\in \nc}$ of seminorms. Let $A=(a_{j,k})_{j,k\in \nc}$ be an infinite matrix with $0\leq a_{j,k}\leq a_{j,k+1}$ and $\displaystyle \sup_{k\in \nc} a_{j,k}>0$ for all $j,k\in \nc$. For $1\leq p<\infty$, the Köthe space $\lambda^{p}(A)$ is defined by
			$$\displaystyle\lambda^{p}(A):=\bigg\lbrace x=(x_{n})_{n\in \nc} :\|x\|_{k}=\bigg(\sum_{j\in \mathbb{N}}|x_{j}a_{j,k}|^{p}\bigg)^{\frac{1}{p}}<+\infty \;\;\text{for all}\;\; k\in \nc\bigg\rbrace,$$
			and for $p=\infty$, the Köthe space $\lambda^{\infty}(A)$ is defined by
			$$\displaystyle \lambda^{\infty}(A):=\bigg\lbrace x=(x_{n})_{n\in \nc} :\|x\|_{k}=\sup_{j\in \nc}|x_{j}|a_{j,k}<+\infty \quad\text{for all}\;\; k\in \nc\bigg\rbrace.$$
			Equipped with the seminorms $\|\cdot\|_{k}$, these spaces are Fr\'echet spaces. If the space $\lambda^p(A)$ is nuclear for some $p$, $1\leq p\leq \infty$, then $\lambda^p(A)$ is nuclear for all $p$, $1\leq p\leq \infty$, and the spaces $\lambda^p(A)$ corresponding to different values of p are the same, see \cite[Proposition 28.16]{vogt4}.
			
			Let $\alpha=(\alpha_{n})_{n\in \nc}$ be an increasing sequence tending to infinity. A power series space of finite type is defined by
			$$\Lambda_{0}\left(\alpha\right):=\left\{x=\left(x_{n}\right)_{n\in \mathbb{N}}: \;\left\|x\right\|_{k}=\sum^{\infty}_{n=1}\left|x_{n}\right|e^{-{1\over k}\alpha_{n}}<\infty \textnormal{ for all } k\in \mathbb{N}\right\},$$
			and a power series space of infinite type is defined by
			$$\displaystyle \Lambda_{\infty}\left(\alpha\right):=\left\{x=\left(x_{n}\right)_{n\in \mathbb{N}}:\; \left\|x\right\|_{k}=\sum^{\infty}_{n=1}\left|x_{n}\right|e^{k\alpha_{n}}<\infty \textnormal{ for all } k\in \mathbb{N}\right\}.$$
			Power series spaces constitute a significant class of Köthe spaces, encompassing the spaces of holomorphic functions $H(\mathbb{D}^{d})$ on the polydisc or the spaces of entire functions $H(\mathbb{C}^{d})$.  The sequence $\alpha$ is called stable if $$\sup_{j\in \nc} \frac{\alpha_{j+1}}{\alpha_j}<\infty .$$
			%%%%%%%%%%%%%%%%%%%%%%%%%%%%%%%%%%%%%%%%%%%%%%%%%%%%%%%%%%%%%%%%%%%%%%%%%%%%%%%%%%%%%%%%%%%%%%%%%%%%%%%%%%%%%%%%%%%%%%%%%%%%%%%%%%%%
			%%%%%%%%%%%%%%%%%%%%%%%%%%%%%%%%%%%%%%%%%%%%%%%%%%%%%%%%%%%%%%%%%%%%%%%%%%%%%%%%%%%%%%%%%%%%%%%%%%%%%%%%%%%%%%%%%%%%%%%%%%%%%%%%%%%%%%%%%%%%%%%%%%%%%%%%
			
			%\subsection{Operators}
			A linear operator $T:E\to F$ is called continuous if for every $k\in \nc$, there exist an $m_{k}\in \nc$ and a $C_{k}>0$ such that
			$$\|T x\|_{k}\leq C_{k}\|x\|_{m_{k}}$$
			for every $x\in E$. The space of linear continuous operators is denoted by $L(E,F)$. A linear map $T: E\to F$ is called bounded, respectively compact, if $T(U)$ is bounded, respectively precompact,
			in F, where $U$ is a neighborhood of zero in E.
			
			A quasi-diagonal operator $T:\lambda^{p}(A)\to \lambda^{q}(B)$ is an operator of the form $Te_{n}=t_{n}e_{\sigma(n)}$ for every $n\in \nc$, where $(t_{n})_{n\in \nc}$ is a sequence of scalars and $\sigma:\nc\to \nc$ is an injective map, here the symbol $e_{n}$ refers to the sequence $(0,0,\dots,0,1,0,\dots)$ where the n$^{th}$ term is 1, and all other terms are 0.
			%%%%%%%%%%%%%%%%%%%%%%%%%%%%%%%%%%%%%%%%%%%%%%%%%%%%%%%%%%%%%%%%%%%%%%%%%%%%%%%%%%%%%%%%%%%%%%%%%%%%%%%%%%%%%%%%%%%%%%%%%%
			
			The characteristic of  continuity map of a linear operator $T$, denoted by $\pi_{T}$, is defined as:
			$$\pi_T(k):= \inf\lbrace r\in \nc: \sup_{\|x\|_r\leq 1}\|Tx\|_k<\infty \rbrace= \inf_{r \in \nc}\|T\|_{k,r}$$
			where $\displaystyle \|T\|_{k,r}:= \sup_{\|x\|_r\leq 1}\|Tx\|_k.$ We want to note that an operator $T$ is bounded if and only if 
			its characteristic of continuity is bounded. In general, it is not possible to put an estimate to the map $\pi_T$. The tameness of an operator provides control over the unpredictable behavior of its characteristic of continuity  is defined as follows:

			%%%%%%%%%%%%%%%%%%%%%%%%%%%%%%%%%%%%%%%%%%%%%%%%%%%%%%%%%%%%%%%%%%%%%%%%%%%%%%%%%%%%%%%%%%%%%%%%%%%%%%%%%%%%%%%%%%%%%%
			
			\begin{defin}
				Let $S:\nc\to\nc$ be a non-decreasing function and  $T\in L(E,F)$, where $E$ and $F$ are graded Fr\'echet spaces. We say that $T$ is \textbf{$S$-tame} if there exists a $k_0\in\nc$ such that  $$\hspace{1.9in}\pi_T(k)\leq S(k)\hspace{1.4in}\forall k\geq k_{0}.$$
				We will call operator $T$ \textbf{tame} if there exists a function S such that T is $S$-tame. If $S$ is a linear function, then $T$ is referred to as \textbf{linearly tame}. 
			\end{defin}
			Vogt showed that every continuous linear map between finite type power series spaces is linearly tame (\cite{vogt5}). In \cite{D1}, Do\u{g}an established the conditions under which all Toeplitz operators defined between power series spaces are $S$-tame. In \cite{D2}, Do\u{g}an demonstrated that all Hankel operators defined between power series spaces are compact and concluded that these operators are always tame.
			%%%%%%%%%%%%%%%%%%%%%%%%%%%%%%%%%%%%%%%%%%%%%%%%%%%%%%%%%%%%%%%%%%%%%%%%%%%%%%%%%%%%%%%%%%%%%%%%%%%%%%%%%%%%%%%%%%%%
			
			\begin{defin} Let $E$ and $F$ be Fr\'echet spaces. The pair $(E,F)$ is called \textbf{tame}, and denoted by $(E,F)\in \mathfrak{T}$ if there exists a non-decreasing function $S:\mathbb{N}\to \mathbb{N}$ such that for every continuous linear operator $T:E\to F$
				$$\pi _T(k)\leq S(k)$$ 
				for sufficiently large k. Equivalently, $(E,F)$ is tame if there exists a sequence of increasing functions $S_\alpha:\nc\to \nc$, $\alpha\in \nc$ such that for every continuous linear operator $T:E\to F$, there exists an $\alpha $ satisfying
				$$\pi_T(k)\leq S_\alpha(k)$$
				for sufficiently large k. 
				If $E=F$, then $E$ is referred to as a \textbf{tame space}.
				Additionally, the pair $(E,F)$ is called \textbf{linearly tame} if $S_\alpha \nearrow S$, where $S$ is a linear function.
			\end{defin}
			
			%{\color{red} For example, when we take $E=F=\Lambda_1(\beta)$, then we can take $S_\alpha(k)=\alpha k$. {\color{red} A proof for this result can be found in ~\cite{vogt3}.}}
			%%%%%%%%%%%%%%%%%%%%%%%%%%%%%%%%%%%%%%%%%%%%%%%%%%%%%%%%%%%%%%%%%%%%%%%%%%%%%%%%%%%%%%%%%%%%%%%%%%%%%%%%%%%%%%%%%%%%%%%%%%%%%%%%%%%%%%%%%5
			
			In 1.1 Proposition of \cite{vogt1}, Dubinsky and Vogt gave a characterization of tame pairs in Fr\'echet spaces as follows:
			
			\begin{pro}\label{sec:tamedef}
				Let $E$ and $F$ be Fr\'echet spaces and $S_{\alpha}:\nc\to \nc$ be nondecreasing functions for $\alpha=1,2,\dots$. The following are equivalent:
				\begin{enumerate}
					\item $(E,F)$ is a tame pair; that is, $\displaystyle \bigcup_{\alpha}L_{\alpha}(E,F)=L(E,F)$ where $L_{\alpha}(E,F)=\big\lbrace T\in L(E,F): \|T\|_{k,S_{\alpha}
						(
						k)}<+\infty \;\; \text{for all}\;\; k\in \nc \big\rbrace$.
					\item For every sequence $K(n)$ of positive integers tending to infinity, there exists an $\alpha$ such that for every $k$, we have $n_0$ and $C$ with 
					$$\|T\|_{k,S_\alpha(k)}\leq C\sup_{n\leq n_0}\|T \|_{n,K(n)}$$
					for every $T\in L(E,F).$
				\end{enumerate}
			\end{pro}
			%%%%%%%%%%%%%%%%%%%%%%%%%%%%%%%%%%%%%%%%%%%%%%%%%%%%%%%%%%%%%%%%%%%%%%%%%%%%%%%%%%%%%%%%%%%%%%%%%%%%%%%%%%%%%%%%%%%%%%%%%%%%%%%%%%%%%%%%5
			
			In Theorem 2.5 of \cite{pisek2}, Piszczek characterized tame pairs in K\"othe spaces as follows:
			
			\begin{theorem}\label{thm1}	Let $A=(a_{j,n})_{j,n\in \nc}$ and $B=(b_{j,n})_{j,n\in \nc}$ be K\"othe matrices. The following are equivalent:
				\begin{enumerate}
					\item $(\lambda^p(A),\lambda^q(B))\in \mathfrak{T}$ for every pair $(p,q)$, $1\leq p,q\leq \infty$, where $p=1$ or $q=\infty$;
					
					\item $(\lambda^p(A),\lambda^q(B))\in \mathfrak{T}$ for some pair $(p,q)$, $1\leq p,q\leq \infty$, where $p=1$ or $q=\infty$;

					\item there exists an increasing function $\psi:\nc\to\nc$, tending to infinity such that for any other increasing function $\phi :\mathbb{N}\to \mathbb{N}$ tending to infinity we have:
					\begin{equation*}\label{E1}
						\exists k_{0}\;\; \forall m\geq k_{0} \;\; \exists n, C_{m} \;\; \forall k,j\in \mathbb{N}:\hspace{0.2in}
						\frac{b_{j,m}}{a_{k,\psi(m)}} \leq C_{m} \max_{1\leq p\leq n} \frac{b_{j,p}}{a_{k,\phi(p)}}.
					\end{equation*}
				\end{enumerate}	
			\end{theorem}
			
			The third statement of Theorem 2.3 can also be expressed as follows:
			\begin{itemize}
				\item[$3^{\prime}$.] There exists a sequence of  non-decreasing functions $S_\alpha:\nc\to\nc$, $\alpha\in \nc$, such that for each $\phi :\mathbb{N}\to \mathbb{N}$ there exists an $\alpha \in \mathbb{N} $ so that for each $m\in \mathbb{N}$ there exist an $n\in \mathbb{N}$ and a $C>0$ with 	
				$$\frac{b_{j,m}}{a_{k,S_{\alpha}(m)}} \leq C_{m} \max_{1\leq p\leq n} \frac{b_{j,p}}{a_{k,\phi(p)}}$$
				holds for every  $k,j\in \nc$. 
			\end{itemize}
			
			%%%%%%%%%%%%%%%%%%%%%%%%%%%%%%%%%%%%%%%%%%%%%%%%%%%%%%%%%%%%%%%%%%%%%%%%%%%%%%%%%%%%%%%%%%%%%%%%%%%%%%%%%%%%%%%%%%%%%%%%%%%%%%%%%%%
			\section{Characterization of tame Fr\'echet spaces using quasi-diagonal operators}

			Dragilev \cite{drag} and Nurlu \cite{nurlu} proved that the existence of a continuous linear unbounded operator between nuclear K\"othe spaces implies the existence of a continuous unbounded quasi-diagonal operator between them.
			By using Theorem \ref{thm1}, we can state a similar result about the tameness of two K\"othe spaces. In the proof of the Theorem \ref{sec:qdoperator}, we closely follow the proof of Corollary 2 in \cite{djakov1}.
			
			\begin{theorem}\label{sec:qdoperator}
				$(\lambda^p(A),\lambda^q(B))$ is tame for every pair $(p,q)$, $1\leq p,q\leq \infty$, where $p=1$ or $q=\infty$ if and only if every continuous quasi-diagonal operator from $\lambda^p(A)$ to $\lambda^q(B)$ is tame. 
			\end{theorem}
			
			\begin{proof}
				It is obvious that if $(\lambda^p(A),\lambda^q(B))$ is tame, then each continuous quasi-diagonal operator from $\lambda^p(A)$ to $\lambda^q(B)$ is tame. So, it is enough to prove that if $(\lambda^p(A),\lambda^q(B))$ is not tame, there exists a quasi-diagonal operator $T:\lambda^{p}(A)\to \lambda^{q}(B)$ which is not tame. Let assume that $(\lambda^p(A),\lambda^q(B))$ is not tame. Then condition 3 in Theorem \ref{thm1} fails. This means that for any increasing function $\psi:\nc\to \nc$, there exists a map $\phi:\nc\to \nc$ such that for all $k\in \nc$ there exists $m_{k}\in \nc$ so that for all $n\in \mathbb{N}$ there are some $i_n, \nu_n\in \nc$ satisfying
				\begin{equation}\label{E2}\frac{b_{\nu_n, m_{k}}}{a_{i_n, \psi(m_{k})}}\geq n \max_{1\leq q\leq n} \frac{b_{\nu_n, q}}{a_{i_n, \phi(q)}}.
				\end{equation}
				Without loss of generality, we can assume the sequences $(i_n)_{n\in \nc}$ and $(\nu_n)_{n\in \nc}$ are strictly increasing, since if they are not, we can find strictly increasing subsequences and work with them. 
			
				Let us take
				$$\hspace{1.4in}t_n^{-1}:= \max_{1\leq q\leq n}\frac{b_{\nu_n, q}}{a_{i_n, \phi(q)}}\hspace{1.5in}\forall n\in \nc$$
				and define a quasi-diagonal operator $T:\lambda^{p}(A)\to \lambda^{q}(B)$ by:	
				\begin{displaymath}
					Te_i= \left\{ \begin{array}{ll}
						t_ne_{\nu_n} & \textrm{if $i=i_n$}\\
						0 & \textrm{if $i\neq i_n$}
					\end{array} \right.
				\end{displaymath}
				It is easy to see that $T$ is continuous. On the other hand, (\ref{E2}) says that for each $k$ there exists $m_{k}$ such that for all $n\in \nc$ we have
				\begin{equation*}
					\frac{\|Te_{i_n}\|_{m_{k}}}{\|e_{i_n}\|_{\psi(m_{k})}} = 
					\frac{1}{\displaystyle \max_{1\leq q\leq n}\frac{b_{\nu_n, q}}{a_{i_n,\phi(q)}}}\cdot \frac{b_{\nu_n,m_{k}}}{a_{i_n,\psi(m_{k})}}  \geq n
				\end{equation*}
				Since $\psi$ is an arbitrary function, this means that $T$ is not tame. This completes the proof. 
			\end{proof}
			
			The proof of the Theorem \ref{sec:qdoperator} characterizing tameness can be adapted by taking $\psi$ linear to establish linear tameness. Therefore we have:
			
			\begin{cor}\label{C1}
				$(\lambda^p(A),\lambda^q(B))$ is linearly tame for every pair $(p,q)$, $1\leq p,q\leq \infty$, where $p=1$ or $q=\infty$ if and only if every continuous quasi-diagonal operator from $\lambda^p(A)$ to $\lambda^q(B)$ is linearly tame. 
			\end{cor}
			%%%%%%%%%%%%%%%%%%%%%%%%%%%%%%%%%%%%%%%%%%%%%%%%%%%%%%%%%%%%%%%%%%%%%%%%%%%%%%%%%%%%%%%%%%%%%%
			
			\section{Tameness of Pairs of Power Series Spaces}
			In this section, we will thoroughly examine the tameness of pairs of power series spaces:
			Zaharjuta's work in \cite{Z}  demonstrated that every continuous linear operator from $\Lambda_{0}(\alpha)$ to $\Lambda_{\infty}(\beta)$ is compact; consequently, in Proposition \ref{P1} we can state that the pair
			$(\Lambda_{0}(\alpha), \Lambda_{\infty}(\beta))$ is tame.  Although Vogt had previously shown the tameness of the pair $(\Lambda_{0}(\alpha), \Lambda_{0}(\beta))$, we present an alternative proof in Theorem \ref{T2}. Nyberg characterized the tameness of the pair $(\Lambda_{\infty}(\alpha), \Lambda_{0}(\beta))$ as in Theorem \ref{T3}, and we demonstrate that the pair $(\Lambda_{\infty}(\alpha), \Lambda_{0}(\beta))$  is not tame if either $\alpha$ or $\beta$ is stable in Theorem \ref{T4}. Furthermore, we provide the proof of Theorem \ref{T1}, which characterizes the tameness of the pair 
			and $(\Lambda_{\infty}(\alpha), \Lambda_{\infty}(\beta))$ was stated without proof by Nyberg. Finally, Theorem \ref{T6} investigates the tameness of the Cartesian products of power series spaces. We can summarize all these results as in Table 1.

			\begin{center}
				\scalebox{0.85}{
					\begin{tabular}{ | c | c | c | c | c | } 
						\hline
						\multirow{2}{6em}{} & $\Lambda_0(\beta)$ & $\Lambda_0(\beta)$ & $\Lambda_\infty(\beta)$ & $\Lambda_\infty(\beta)$ \\
						& stable& non-stable & stable & non-stable \\
						\hline 
						\multirow{2}{6em}{\centering{$\Lambda_0(\alpha)$} \\ stable} & Tame  & Tame& Tame& Tame \\ 
						&Theorem \ref{T2} & Theorem \ref{T2} & Proposition \ref{P1} & Proposition \ref{P1}\\
						\hline
						\multirow{2}{6em}{\centering{$\Lambda_{0}(\alpha)$}\\ non-stable} & Tame & Tame & Tame & Tame \\ 
						&Theorem \ref{T2}   & Theorem \ref{T2}  & Proposition \ref{P1}& Proposition \ref{P1}\\
						\hline
						\multirow{2}{6em}{\centering{$\Lambda_\infty(\alpha)$} \\ stable} & Non-Tame  & Non-Tame & Tame iff Bounded& Tame iff Bounded \\ 
						& Theorem \ref{T4} & Theorem \ref{T4} & Proposition \ref{P2}& Proposition \ref{P2} \\
						\hline
						\multirow{2}{6em}{\centering{$\Lambda_\infty(\alpha)$} \\ non-stable} & Non-Tame & Tame & Tame iff Bounded & Tame\\ 
						& Theorem \ref{T4}  & Theorem \ref{T3} & Proposition \ref{P2} &   Theorem \ref{T1}\\
						\hline
				\end{tabular}}
				~\\~\\
				\textbf{Table 1}
			\end{center}
			
			%%%%%%%%%%%%%%%%%%%%%%%%%%%%%%%%%%%%%%%%%%%%%%%%%%%%%%%%%%%%%%%%%%%%%%%%%%%%%%%%%%%%%%%%%%%%%%%%%%%%%%%%%%%%%%%%%%%%%%%%%%%
			Before presenting the main results, let us state the following lemma, which will be used in the proofs of  Theorem \ref{T2} and Theorem \ref{sec:composition}, and the operators constructed within its proof will be applied in the subsequent sections of the article.
			
			\begin{lemma}\label{sec:operator} For every $T\in L(\Lambda_r(\alpha),\Lambda_r(\beta))$, $r\in  \lbrace 0,\infty\rbrace$, there exist an increasing sequence $\gamma=(\gamma_n)_{n\in\nc}$ tending to infinity, and an operator $R\in L(\Lambda_r(\gamma),\Lambda_r(\gamma))$ such that the range of $T$ is isomorphic to the range of $R$. 
			\end{lemma}
			\begin{proof}  Since $\alpha$ and $\beta$ are increasing sequences, one can construct a new increasing sequence $\gamma$ such that no term of $\alpha$ or $\beta$ is omitted in $\gamma$. Then, there exist two increasing sequences $(t_n)_{n\in \mathbb{N}}, (s_n)_{n\in \mathbb{N}}\subseteq \mathbb{N}$ such that 
				$$\gamma_{t_n}=\alpha_n, \hspace{0.3in}\text{and}\hspace{0.3in} \gamma_{s_n}=\beta_n$$
				for all $n\in \nc$, and $\lbrace (\gamma_n)_{n\in\nc}\rbrace=\lbrace (\gamma_{t_n})_{n\in\nc}\rbrace \cup \lbrace (\gamma_{s_n})_{n\in\nc}\rbrace$. Let us define an operator $T_1:\Lambda_r(\gamma)\to \Lambda_r(\alpha) $ by $$T_1((x_n)_{n\in\nc}):=(x_{t_n})_{n\in\nc}.$$ Then $T_{1}$ is a continuous operator since we have
				\begin{equation}\label{E4}
					\|T_1x\|_k^2= \sum_{n\in\nc}|x_{t_n}|^2e^{2r_k\alpha_n}= \sum_{n\in\nc}|x_{t_n}|^2e^{2r_k\gamma_{t_n}}\leq \sum_{n\in\nc}|x_n|^2e^{2r_k\gamma_n} =\|x\|_k^2
				\end{equation}
				for every $k\in \mathbb{N}$, here  $r_{k}$'s are taken as $-\frac{1}{k}$ for $r=0$ and as $k$ for $r=\infty$. Furthermore, $T_{1}$ is surjective. Indeed, for every $x=(x_n)_{n\in\nc}\in \Lambda_r(\alpha)$, we define $u=(u_i)_{i\in\nc}\in \Lambda_{r}(\gamma)$ by
				\begin{equation*}
					u_i = \left\{
					\begin{array}{cc}
						x_n & \text{if } i=t_n,\\
						0 & \text{otherwise}. 
					\end{array} \right.
				\end{equation*}
				and it is clear that $T_1u=x$. Let us define an operator $T_2:\Lambda_r(\beta)\to \Lambda_r(\gamma) $ by $$T_2((y_n)_{n\in\nc}):=(z_i)_{i\in\nc}$$ where
				\begin{equation*}
					z_i = \left\{
					\begin{array}{cc}
						y_n & \text{if } i=s_n,\\
						0 & \text{otherwise}. 
					\end{array} \right.
				\end{equation*}
				Then $T_{2}$ is a continuous operator since we have
				\begin{equation}\label{E3}
					\|T_2y\|_k^2= \sum_{i\in\nc}|z_i|^2e^{2r_k\gamma_i}
					= \sum_{n\in\nc}|y_n|^2e^{2r_k\gamma_{s_n}}\leq \sum_{n\in\nc}|y_n|^2e^{2r_k\beta_n} =\|y\|_k^2
				\end{equation}
				for all $k\in \nc$. $T_{2}$ is injective as its kernel contains only the zero.  
				Now we define the operator 
				$R\in L(\Lambda_r(\gamma),\Lambda_r(\gamma))$ as
				$$R:= T_2\circ T \circ T_1.$$ 
				Since $T_{1}$ is surjective and $T_{2}$ is injective continuous operators, the range of $R$ is isomorphic to the range of $T$. 
			\end{proof}
			
			%%%%%%%%%%%%%%%%%%%%%%%%%%%%%%%%%%%%%%%%%%%%%%%%%%%%%%%%%%%%%%%%%%%%%%%%%%%%%%%%%%%%%%%%%%%%%%%%%%
			Zaharjuta in \cite{Z} showed that every continuous linear operator $T:\Lambda_0(\alpha)\to \Lambda_\infty(\beta)$ is compact, as a direct consequence we have the following proposition:
			\begin{pro}\label{P1} For every sequence $\alpha$ and $\beta$, the pair $(\Lambda_{0}(\alpha),\Lambda_\infty(\beta))$ is tame.
			\end{pro}
			%%%%%%%%%%%%%%%%%%%%%%%%%%%%%%%%%%%%%%%%%%%%%%%%%%%%%%%%%%%%%%%%%%%%%%%%%%%%%%%%%%%%%%%%%%%%%%%%%%
			%~\\
			%{\color{gray} Piszczek Dubinsky ve Vogt'a atıf vermiş, Nyberg ve Mitiagin'den bahsetmiş.
				%~\\
				%Nyberg; 1.2. Proposition: Pairs of nuclear finite type power series spaces are tame.
				%~\\ 
				%Piszczek; Theorem 4.1: Let $\alpha$ be stable and let X be a Fr\'echet space. Assume that either X is nuclear or $X =\lambda^{1}(A)$ or $p = +\infty$. The following are equivalent:
				%\begin{itemize}
				%\item[(i)] a pair $(X,\Lambda^p_0(\alpha))$ is tame,
				%\item[(ii)] ...
				%\item[(iii)] $X \in (\overline{\Omega})$.
				%\end{itemize}
				%
				%~\\~\\
				%Vogt Operators Between Fr\'echet spaces 2.3 Theorem: Every continuous linear map from $\Lambda_{0}(\alpha)$ to $\Lambda_{0}(\beta)$ is linearly tame. 
				%~\\~\\
				%Aşağıdaki teoremde ise nükleerlik yok, diziler üzerinde stable olma koşulu yok ve ispatta quasi-diagonal operatörler kullanılıyor.
				%}
			%{\color{red} Buraya bir kaç cümle gelecek!}
			
			Vogt demonstrated in Lemma 5.1 of \cite{vogt3} that 
			every finite type power series space is tame. Vogt also showed that the pair $(\Lambda_0(\alpha),\Lambda_0(\beta))$ is tame (and linearly tame) in 2.2 Theorem (and 2.3 Theorem) of \cite{vogt4}. Below, it is shown that the pair $(\Lambda_0(\alpha),\Lambda_0(\beta))$ is tame by employing a different technique.
			
			\begin{theorem}\label{T2} The pair
				$(\Lambda_0(\alpha),\Lambda_0(\beta))$ is always tame.
			\end{theorem}
			\begin{proof} Let us assume that a pair $(\Lambda_0(\alpha),\Lambda_0(\beta))$ is not tame and $S:\nc\to \nc$ be an increasing function. Then by Theorem \ref{sec:qdoperator}, there exists a quasi-diagonal operator $T\in L(\Lambda_0(\alpha),\Lambda_0(\beta))$ which is not $S$-tame. This means that  for every $k,j \in\nc$, there exist  $m_{k}\geq k$ and $x_{j}=(x_{j,n})_{n\in\nc} \in \Lambda_{0}(\alpha)$ such that 	$$\frac{\|Tx_j\|_{m_k}}{\|x_j\|_{S(m_k)}}\geq j.$$
				We construct the sequence $\gamma$ and the operator $R\in L(\Lambda_0(\gamma),\Lambda_0(\gamma))$ as in Lemma \ref{sec:operator}. We define the sequence $u_j=(u_{j,i})_{i\in\nc}\in \Lambda_{0}(\gamma)$ as
				\begin{equation*}
					u_{j,i} = \left\{
					\begin{array}{cc}
						x_{j,n} & \text{if } i=t_n,\\
						0 & \text{otherwise }. 
					\end{array} \right.
				\end{equation*}
				Then $T_{1}(u_{j})=x_{j}$ for every $j\in \nc$ and we have
				\begin{equation}\label{E5}
					\begin{split}
						\|Ru_j\|^2_k&=\sum_{n\in\nc}|(Ru_j)_{n}|^2e^{-2k\gamma_n}=\sum_{n\in\nc}|(T_2\circ T(T_1u_j))_n|^2e^{-2k\gamma_n}\\
						&=\sum_{n\in\nc}|(T_2(Tx_j))_n|^2e^{-2k\gamma_n}= \sum_{n\in\nc}|(Tx_j)_n|^2e^{-2k\gamma_{s_n}}\\ &=\sum_{n\in\nc}|(Tx_j)_n|^2e^{-2k\beta_n}=\|Tx_j\|^2_k.
					\end{split}
				\end{equation}
				Similarly, we can write
				\begin{equation}\label{E6}
					\begin{split}
						\|u_j\|^2_{S(k)}&=\sum_{i\in\nc}|u_{j,i}|^2e^{-2S(k)\gamma_i}=\sum_{n\in\nc}|x_{j,n}|^2e^{-2S(k)\gamma_{t_n}}\\
						&=\sum_{n\in\nc}|x_{j,n}|^2e^{-2S(k)\alpha_n}=\|x_j\|^2_{S(k)}.
					\end{split}
				\end{equation}
				So we have 
				$$\frac{\|Ru_j\|_{m_k}}{\|u_j\|_{S_\alpha(m_k)}}\geq j.$$
				Therefore $R$ is a not $S$-tame operator on $\Lambda_0(\gamma)$. Since $S$ was arbitrary, we conclude that $\Lambda_0(\gamma)$ is not tame, this is a contradiction. Therefore, every pair $(\Lambda_{0}(\alpha), \Lambda_{0}(\beta))$ is tame.
			\end{proof}
			
			Since the proof above also holds for linear 
			$S$ functions, we can state the following result.
			\begin{cor} The pair $(\Lambda_0(\alpha),\Lambda_0(\beta))$ is linearly tame.
			\end{cor}
			%%%%%%%%%%%%%%%%%%%%%%%%%%%%%%%%%%%%%%%%%%%%%%%%%%%%%%%%%%%%%%%%%%%%%%%%%%%%%%%%%%%%%%%%%%%%%%%%%%%%%%%%%%%%%%%%%%%%%%%%%%%%%%%%%%%%%%%%%%%%%%%
			
			The following theorem is proved by Nyberg 
			and stated as 1.7 Theorem in  \cite{nyberg}:
			\begin{theorem}\label{T3}
				The following conditions are equivalent.
				\begin{itemize}
					\item[(i)] Every operator from $\Lambda_\infty(\alpha)$ to $\Lambda_0(\beta)$ is compact.
					\item[(ii)] $(\Lambda_\infty(\alpha),\Lambda_0(\beta))$ is tame.
					\item[(iii)] $(\Lambda_\infty(\alpha),\Lambda_0(\beta))$ is linearly tame.
					\item[(iv)] $M_{\beta\alpha}$; the set of finite limit points of $(\beta_i/\alpha_j)_{i,j\in\nc}$ is bounded.
				\end{itemize}
			\end{theorem}
			
			\begin{rem}\label{rem1} We want to note that if both $\alpha$ and $\beta$ are non-stable, the set of finite limit points of  $(\beta_i/\alpha_j)_{i,j\in\nc}$ is bounded.
			\end{rem}
			\begin{theorem}
					If $\alpha$ and $\beta$ is non-stable then $(\Lambda_\infty(\alpha),\Lambda_0(\beta))$ is tame.
			\end{theorem}
			%%%%%%%%%%%%%%%%%%%%%%%%%%%%%%%%%%%%%%%%%%%%%%%%%%%%%%%%%%%%%%%%%%%%%%%%%%%%%%%%%%%%%%%%%%%%%%%%%%%%%%
			
			We can draw further conclusions when $\alpha$ or $\beta$ is stable.
			
			\begin{theorem}\label{T4}
				If $\alpha$ or $\beta$ is stable then $(\Lambda_\infty(\alpha),\Lambda_0(\beta))$ is not tame.
			\end{theorem}
			
			\begin{proof} Let us assume that $\alpha$ is stable and $(\Lambda_\infty(\alpha),\Lambda_0(\beta))$ is tame. Theorem 4.10 of \cite{pisek2} gives us that every continuous linear operator from $\Lambda_{\infty}(\alpha)$ to $\Lambda_{0}(\beta)$ is bounded. By using Satz 3.2 of \cite{vogt2}, we can say that $\Lambda_0(\beta)$ has the property $(LB_\infty)$. However, the remark given in \cite{vogt2} on page 190 indicates that the property $(LB_{\infty})$ is stronger than the property $(DN)$. In this case, $\Lambda_0(\beta)$ has the property $(DN)$, which leads to a contradiction.
				
				Let us consider the case that $\beta$ is stable and  assume that  $(\Lambda_\infty(\alpha),\Lambda_0(\beta))$ is tame. Theorem 4.1 of \cite{pisek2} says that $\Lambda_\infty(\alpha)$ has the property $(\overline{\Omega})$, this is a contradiction. Therefore, we can say that $(\Lambda_\infty(\alpha),\Lambda_0(\beta))$ is not tame if $\alpha$ or $\beta$ is stable.
			\end{proof}
			%%%%%%%%%%%%%%%%%%%%%%%%%%%%%%%%%%%%%%%%%%%%%%%%%%%%%%%%%%%%%%%%%%%%%%%%%%%%%%%%%%%%%%%%%%%%%%%%%%%%%%%%%%%%%%%%%%%%%%%%%%%%%%%%%%%%%%%%%%%%%%%
			
			The tameness of power series spaces of infinite type was investigated by Dubinsky and Vogt in \cite[1.3 Theorem]{vogt1}. Dubinsky and Vogt proved in Theorem 1.3 that the tameness of $\Lambda_{\infty}(\alpha)$, the linearly tameness of $\Lambda_{\infty}(\alpha)$ and the boundedness of the set of finite limit points of $(\alpha_i/\alpha_j)_{i,j\in\nc}$ are equivalent.
			The following theorem was presented without a proof in \cite[1.4 Proposition]{nyberg} and our objective is to establish its proof.
			
			\begin{theorem}\label{T1}
				The following statements are equivalent:
				\begin{itemize}
					\item[(i)] $(\Lambda_\infty(\alpha),\Lambda_\infty(\beta))$ is tame.
					\item[(ii)] $(\Lambda_\infty(\alpha),\Lambda_\infty(\beta))$ is linearly tame.
					\item[(iii)] $M_{\beta\alpha}$; the set of finite limit points of $(\beta_i/\alpha_j)_{i,j\in\nc}$ is bounded. 
				\end{itemize}
			\end{theorem}

			%%%%%%%%%%%%%%%%%%%%%%%%%%%%%%%%%%%%%%%%%%%%%%%%%%%%%%%%%%%%%%%%%%%%%%%%%%%%%%%%%%%%%%%%%%%%%%%%%%%%%%%%%%%%%%%%%%%%%%%%%%%%%%%%%
			First, let us consider the case where $(i)\Leftrightarrow (ii)$:
			
			\begin{pro}\label{PT1}
				$(\Lambda_\infty(\alpha),\Lambda_\infty(\beta))$ is tame if and only if $(\Lambda_\infty(\alpha),\Lambda_\infty(\beta))$ is linearly tame.
			\end{pro}	
			\begin{proof} It is clear that if $(\Lambda_\infty(\alpha),\Lambda_\infty(\beta))$ is linearly tame, then $(\Lambda_\infty(\alpha),\Lambda_\infty(\beta))$ is tame. Let us assume that $(\Lambda_\infty(\alpha),\Lambda_\infty(\beta))$ is tame, but not linearly tame. Then there exists a tame operator $T\in L(\Lambda_\infty(\alpha),\Lambda_\infty(\beta))$ which is not linearly tame, that is, there exists an increasing function $S:\nc \to \nc$ and $C>0$ such that for every sufficiently large $k$, there exists a $C_k>0$ with 
				$$\|Tx\|_k\leq C \|x\|_{S(k)}$$
				for all $x\in \Lambda_\infty(\alpha)$, but for any $a,b\in\mathbb{R}$, $S_{ab}(k)=ak+b$, and  any $k\in\nc$, we can find $x_j=(x_{j,n})_{n\in\nc}$ satisfying
				$$\frac{\|Tx_j\|_k}{\|x_j\|_{S_{ab}(k)}}\geq j$$
				for all  $j\in\nc$. We again consider the sequence $\gamma$ and the operator $R\in L(\Lambda_0(\gamma),\Lambda_0(\gamma))$ defined as in Lemma \ref{sec:operator}. By using (\ref{E4}) and (\ref{E3}), we have
				$$\|Ru\|_k=\|T_2\circ T \circ T_1u\|_k\leq \|T\circ T_1u\|_k\leq C \|T_1 u\|_{S(k)}\leq C\|u\|_{S(k)}$$
				for every $u\in \Lambda_\infty(\gamma)$. So, $R$ is a tame operator in $\Lambda_\infty(\gamma)$. On the other hand, for any $u_j=(u_{j,i})_{i\in\nc}$ given by
				\begin{equation*}
					u_{j,i} = \left\{
					\begin{array}{cc}
						x_{j,n} & \text{if } i=t_n,\\
						0 & \text{otherwise},
					\end{array} \right.
				\end{equation*}
				it is easy to demonstrate that 
				$\|Ru_{j}\|_{k}=\|Tx_{j}\|_{k}$ and $\|x_{j}\|_{k}=\|u_{j}\|_{k}$ for every $k\in \nc$ by following  an argument similar to that of equations (\ref{E5}) and (\ref{E6}) and then we have
				$$\|Ru_j\|_{k}=\|Tx_j\|_k\geq j\|x_j\|_{S_{ab}(k)}=j\|u_j\|_{S_{ab}(k)}.$$
				This says that $R$ is not linearly tame. This, in turn, stands in contradiction to Theorem 1.3 of \cite{vogt1}. Therefore $(\Lambda_\infty(\alpha),\Lambda_\infty(\beta))$ is linearly tame.
			\end{proof}
			%%%%%%%%%%%%%%%%%%%%%%%%%%%%%%%%%%%%%%%%%%%%%%%%%%%%%%%%%%%%%%%%%%%%%%%%%%%%%%%%%%%%%%%%%%%%%%%%%%%%%%%%%%%%%%%%%%%%%%%%%%%%%%%%
			
			The following proposition represents the implication $(iii)\Rightarrow (ii)$ in Theorem \ref{T1}.
			
			\begin{pro}\label{sec:flp2}
				If the set of finite limit points of $(\beta_j/\alpha_k)_{j,k\in\nc}$ is bounded, then the pair $(\Lambda_\infty(\alpha),\Lambda_\infty(\beta))$ is linearly tame. 
			\end{pro}
			\begin{proof} Let us assume that the set of finite limit points of $(\beta_j/\alpha_k)_{j,k\in\nc}$ is bounded. It is enough to show that every continuous quasi-diagonal operator is linearly tame by Corollary \ref{C1}. Let $T\in L(\Lambda_\infty(\alpha),\Lambda_\infty(\beta))$ be a quasi-diagonal operator defined by $Te_n=t_ne_{\sigma(n)}$ for every $n\in\nc$.  We want to show that $T$ is linearly tame, that is, there exist $A,B>0$ such that for every sufficiently large $k\in \nc$, there exists $D_k>0$ with 
				$$\frac{\|Te_n\|_k}{\|e_n\|_{Ak+B}}\leq D_k$$
				for all $n\in \nc$. Since $T$ is continuous, for every $k\in \nc$, there exists a $C_{k}>0$ such that
				$$\|Tx\|_{k}\leq C_{k}\|x\|_{\pi_{T}(k)}$$
                for all $x\in \Lambda_{\infty}(\alpha).$
				Let $k>1$ and  $A$ be a strict upper bound for the finite limit points of the set $(\beta_j/\alpha_k)_{j,k\in\nc}$. We consider the intervals $I_1=\lbrack 0, A\rbrack$, $I_2=\lbrack \pi_T(k+1),\infty)\cap\lbrack 2A,\infty)$ and $I_3=\mathbb{R}^+\setminus (I_1\cup I_2)$. 
				
				If $\beta_{\sigma(n)}/\alpha_n\in I_1$, then since $T$ is continuous, there exists a $C_{1}$ such that the inequalities
				\begin{align*}
					C_1&\geq \frac{\|Te_n\|_1}{\|e_n\|_{\pi_T(1)}}=\frac{\|Te_n\|_1}{\|e_n\|_{\pi_T(1)}}\cdot\frac{\|Te_n\|_k}{\|Te_n\|_k}\cdot\frac{\|e_n\|_{Ak+\pi_T(1)}}{\|e_n\|_{Ak+\pi_T(1)}}\\
					&= \frac{\|Te_n\|_k}{\|e_n\|_{Ak+\pi_T(1)}}\cdot e^{(1-k)\beta_{\sigma(n)}}\cdot e^{Ak\alpha_n}\\
					&\geq  \frac{\|Te_n\|_k}{\|e_n\|_{Ak+\pi_T(1)}}\cdot e^{-k\beta_{\sigma(n)}}\cdot e^{Ak\alpha_n}\\
					&\geq \frac{\|Te_n\|_k}{\|e_n\|_{Ak+\pi_T(1)}}\cdot e^{-kA\alpha_n}\cdot e^{Ak\alpha_n}= \frac{\|Te_n\|_k}{\|e_n\|_{Ak+\pi_T(1)}}.
				\end{align*}
				hold. If $\beta_{\sigma(n)}/\alpha_n\in I_2$, then again since $T$ is continuous, there exists a $C_{k+1}$ such that the inequalities
				\begin{align*}
					C_{k+1}&\geq \frac{\|Te_n\|_{k+1}}{\|e_n\|_{\pi_T(k+1)}}=\frac{\|Te_n\|_{k+1}}{\|e_n\|_{\pi_T(k+1)}}\cdot\frac{\|Te_n\|_k}{\|Te_n\|_k}\cdot\frac{\|e_n\|_{Ak+\pi_T(1)}}{\|e_n\|_{Ak+\pi_T(1)}}\\
					&= \frac{\|Te_n\|_k}{\|e_n\|_{Ak+\pi_T(1)}}\cdot e^{\beta_{\sigma(n)}}\cdot e^{(Ak+\pi_T(1)-\pi_T(k+1))\alpha_n}\\
					&\geq  \frac{\|Te_n\|_k}{\|e_n\|_{Ak+\pi_T(1)}}\cdot e^{(Ak+\pi_T(1))\alpha_n}\geq  \frac{\|Te_n\|_k}{\|e_n\|_{Ak+\pi_T(1)}}.
				\end{align*}
				hold.
				Since the set of finite limit points of $(\beta_j/\alpha_k)_{j,k\in\nc}$ is bounded, the set $I_3\cap \lbrace (\beta_j/\alpha_i)\rbrace_{i,j\in\nc}$ is finite, so there are only finitely many n, which we denote by $n_{1}, n_{2}, \dots, n_{m}$, such that $\beta_{\sigma(n)}/\alpha_n\in I_2$ for these n. By defining $D_{k}$ and $B$ as
				$$D_k=\max_{1\leq i\leq m}\lbrace C_{n_i},C_1,C_{k+1}\rbrace, \hspace{0.3in} B=\pi_T(1).$$
				we showed that 
				$$\frac{\|Te_n\|_k}{\|e_n\|_{Ak+B}}\leq D_k$$
				for all $n\in \nc$ . Therefore, $T$ is linearly tame.
			\end{proof}
			%%%%%%%%%%%%%%%%%%%%%%%%%%%%%%%%%%%%%%%%%%%%%%%%%%%%%%%%%%%%%%%%%%%%%%%%%%%%%%%%%%%%%%%%%%%%%%%%%%%%%%%%%%%%%%%%%%%%%%%%%%%%%%%%%%%%%%%%%%%%%%%%
			
			In the proof of Proposition \ref{sec:flp1}, we closely follow and adapt the proof of Theorem 1.6 in \cite{nyberg}.
			\begin{pro}\label{sec:flp1} If $(\Lambda_\infty(\alpha),\Lambda_\infty(\beta))$ is tame, then the set of finite limit points of $(\beta_j/\alpha_k)_{j,k\in\nc}$ is bounded.
			\end{pro}
			
			\begin{proof} Let us assume that $(\Lambda_\infty(\alpha),\Lambda_\infty(\beta))$ is tame, the set of finite limit points of  $(\beta_j/\alpha_k)_{j,k\in\nc}$ is unbounded and $S:\nc\to \nc$ be an increasing function. We may assume that
				$$\hspace{0.7in} S(k+1)-S(k)\geq S(k)-S(k-1)>0\hspace{0.7in}
				\forall\;\; k=2,3,\dots.$$
				By following the same procedure of the proof of Proposition 1.6 in \cite{nyberg}, we can find infinitely many finite limit points $r_k\in \mathbb{R}^+$ with $r_{k+1}>r_k+1$, disjoint infinite sets $M_k\subset \nc$, indices $m(k,n)$, $n\in M_k$, $k\in \nc$, open intervals 
				$$I_k=\bigg( k(k-1)(S(k+1)-S(k)),k(k+1)(S(k+2)-S(k+1))\bigg)$$
				with $\beta_{m(k,n)}/\alpha_n\in I_k$, $n\in M_k$, $k\in \nc$, and
				$$\lim_{\substack{n\in M_k \\ n\to \infty}} \frac{\beta_{m(k,n)}}{\alpha_n}=r_k.$$
				We define the quasi-diagonal operator $T\in L(\Lambda_\infty(\alpha),\Lambda_\infty(\beta))$ as follows:
				\begin{equation*}
					Te_n = \left\{
					\begin{array}{cc}
						e^{S(k+1)\alpha_n}\cdot e^{-k\beta_{m(k,n)}}e_{m(k,n)} &\hspace{0.2in} \text{if } n\in M_k, k\in \nc\\
						0 & \text{otherwise}.
					\end{array} \right.
				\end{equation*}
				Then for every $n\in M_k$ and $k\in \nc$, we have
				$$\|Te_n\|_k=e^{S(k+1)\alpha_n}\cdot e^{-k\beta_{m(k,n)}}\cdot e^{k\beta_{m(k,n)}}=e^{S(k+1)\alpha_n}= \|e_n\|_{S(k+1)}.$$
				Now we prove by induction that for every $l\in \nc$ and for every $n\in \bigcup_{k\in \nc}M_{k}$,
				$$\|Te_{n}\|_{l}\leq \|e_{n}\|_{S(l+1)}.$$
				Let $k\in \nc$ be fixed. Since $\beta_{m(k,n)}\geq k(k-1)(S(k+1)-S(k))\alpha_n$, we can write
				\begin{align*}
					\|Te_n\|_{k-1}&= e^{S(k+1)\alpha_n}\cdot e^{-\beta_{m(k,n)}}\\
					&\leq e^{S(k+1)\alpha_n}\cdot e^{-k(k+1)(S(k+1)-S(k))\alpha_n}\\
					&= e^{\alpha_n(-(k^2-k-1)S(k+1)+(k^2-k)S(k))}\\
					&\leq e^{S(k)\alpha_n}=\|e_n\|_{S(k)}
				\end{align*}
				We have a similar inequality for $k+1$, and by applying induction, we get $\|Te_n\|_j\leq \|e_n\|_{S(j+1)}$ for all $j\in \nc$. So $T$ is a continuous operator and we have	
				\begin{align*}
					\sup_{n\in\nc}\frac{\|Te_n\|_k}{\|e_n\|_{S(k)}}&\geq \sup_{n\in M_k}\frac{\|Te_n\|_k}{\|e_n\|_{S(k)}}\\
					&= \sup_{n\in M_k}\frac{\|e_n\|_{S(k+1)}}{\|e_n\|_{S(k)}}\\
					&=\sup_{n\in M_k}e^{S(k+1)-S(k)\alpha_{n}}\to \infty.
				\end{align*}
				So $T$ is not $S$-tame. Since $S$ was arbitrary, this implies that $(\Lambda_\infty(\alpha),\Lambda_\infty(\beta))$ is not tame.	
			\end{proof}
			
			%{\color{red} Aşağıdaki proof dergi formatına göre düzenlenecek!}
			
			\begin{refproof}{Theorem \ref{T1}}
				(i) $\Leftrightarrow$ (ii) was proven 
				in Proposition \ref{PT1}, 
				(i) $\Rightarrow$ (iii) was proven in Proposition \ref{sec:flp1}, and (iii) $\Rightarrow$ (ii) was proven in Proposition \ref{sec:flp2}. 
			\end{refproof}
			
			%%%%%%%%%%%%%%%%%%%%%%%%%%%%%%%%%%%%%%%%%%%%%%%%%%%%%%%%%%%%%%%%%%%%%%%%%%%%%%%%%%%%%%%%%%%%%%%%%%%%%%%%%%%%%%%%%%%%%%%%%%%%%%%%%%%%%%%%%%%%%%%%%%%%%%%%%%%%%%%%%%%%%%%%%%%%%%%%%%%%%%%%%%%%%%%%
			In Theorems 4.2 and 4.10 in \cite{pisek2}, Piszczek characterized the tameness of the pair $(E, F)$ as the boundedness of all continuous linear operators from $E$ to $F$, provided that one of E or F is an infinite type power series space generated by a stable sequence. The following proposition is a direct consequence of these two theorems.
			
			\begin{pro}\label{P2}
				If $\alpha$ or $\beta$ is stable then $(\Lambda_\infty(\alpha),\Lambda_\infty(\beta))$ is tame if and only if every continuous linear operator between these spaces is bounded.
			\end{pro}
			
			\begin{proof}
				This is a direct consequence of the Theorems 4.2 and 4.10 in \cite{pisek2} applied to $(\Lambda_\infty(\alpha),\Lambda_\infty(\beta))$.
			\end{proof}
			
			By employing Remark \ref{rem1}, we can state:
			\begin{theorem}
					If $\alpha$ and $\beta$ is non-stable then $(\Lambda_\infty(\alpha),\Lambda_\infty(\beta))$ is tame.
			\end{theorem}
			%%%%%%%%%%%%%%%%%%%%%%%%%%%%%%%%%%%%%%%%%%%%%%%%%%%%%%%%%%%%%%%%%%%%%%%%%%%%%%%%%%%%%%%%%%%%%%%%%%%%%
			
			We can apply these results to examine the tameness of cartesian product of power series spaces.
			\begin{theorem}\label{T6} We have the following about the cartesian product of power series spaces:
				\begin{itemize}
					\item 	The cartesian product $\Lambda_0(\alpha)\times \Lambda_0(\beta)$ is always tame.
					\item 	The cartesian product $\Lambda_0(\alpha)\times \Lambda_\infty(\beta)$ is tame if and only if both the sets of finite limit points of $(\beta_i/ \beta_j)_{i,j\in\nc}$ and $(\alpha_i/\beta_j)_{i,j\in\nc}$ are bounded
					\item 	The cartesian product $\Lambda_\infty(\alpha)\times \Lambda_\infty(\beta)$ is tame if and only if the sets of finite limit points of $(\beta_i/\beta_j)_{i,j\in\nc}$, $(\alpha_i/\alpha_j)_{i,j\in\nc}$, $(\beta_i/\alpha_j)_{i,j\in\nc}$ and $(\alpha_i/\beta_j)_{i,j\in\nc}$ are bounded
				\end{itemize}
			\end{theorem}
			
			\begin{proof}
				$T\in L(\Lambda_r(\alpha)\times\Lambda_s(\beta))$ has a matrix representation:
				\begin{equation*}T=
					\begin{bmatrix}
						T_{11} & T_{12} \\
						T_{21} & T_{22}
					\end{bmatrix} : (x,y)\to (T_{11}x+T_{12}y,T_{21}x+T_{22}y),
				\end{equation*}
				where $T_{11}\in L(\Lambda_r(\alpha),\Lambda_r(\alpha))$, $T_{12}\in L(\Lambda_r(\alpha),\Lambda_s(\beta))$, $T_{21}\in L(\Lambda_s(\beta),\Lambda_r(\alpha))$, $T_{22}\in L(\Lambda_s(\beta),\Lambda_s(\beta))$. The result follows.
			\end{proof}
			
			%%%%%%%%%%%%%%%%%%%%%%%%%%%%%%%%%%%%%%%%%%%%%%
\section{Bases of The Range of A Tame Operator Between Power Series Spaces} %{\color{red} of Infinite Type}}
%In this part, we will focus on the range of a tame operator defined between power series spaces of infinite type. 
			Krone \cite{K} introduced a topological condition DS, for nuclear Fréchet spaces $E$ and $F$, which ensures that the range of every continuous linear operator from $E$ to $F$ possesses a basis.  In Theorem 2.1 of \cite{K}, he showed that all continuous linear operators from $\Lambda_{0}(\alpha)$ to $\Lambda_{0}(\beta)$ have ranges with a basis and he also proved that every continuous linear operator from $\Lambda_{\infty}(\alpha)$ to $\Lambda_{\infty}(\beta)$  has a range admitting a basis if and only if the set of finite limit points of $(\beta_{i}/\alpha_{j})_{i,j\in \mathbb{N}}$ is bounded. In the following theorem,  we show that every tame operator from  $\Lambda_{\infty}(\alpha)$ to $\Lambda_{\infty}(\beta)$ has a range admitting a basis, without any restriction on the limit points of $(\beta_{i}/\alpha_{j})_{i,j\in \mathbb{N}}$.		
			%%%%%%%%%%%%%%%%%%%%%%%%%%%%%%%%%%%%%%%%%%%%%%%%%%%%%%%%%%%%%%%%%%%%%%%%%%%%%%%%%%%%%%%%%%%%%%%%%%%%%%%%%%%%%%%%%%%%%%%%%%%%%%%%%%%%%%%%%%%%%%%%%%%

			%{\color{red} Aşağıdaki sonucun sonlu tip uzaylar için hali verilebilir mi? sonlu tip uzaylarda her operatör kompakt ve sonlu tip uzaylarda her complemented subspace yine sonlu tip uzay ama daha başka birşey diyemiyorum. Complemented subspace ile ilgili birşeyler yazılabilir, Vogt Power series space representation}
			
			\begin{theorem}\label{sec:composition}
				The range of every tame operator $T\in L(\Lambda_\infty(\alpha),\Lambda_\infty(\beta))$ has a basis. If $\Lambda_\infty(\beta)$ is nuclear, then the range of $T$ has an absolute basis.
			\end{theorem}
			
			\begin{proof} Let $T\in L(\Lambda_\infty(\alpha),\Lambda_\infty(\beta))$ be a tame operator. In Lemma \ref{sec:operator}, we proved that there exist an increasing sequence $\gamma=(\gamma_{n})_{n\in \nc}$ and operators $T_{1}:\Lambda_{\infty}(\gamma)\to \Lambda_{\infty}( \alpha)$ and $T_{2}:\Lambda_{\infty}(\beta)\to \Lambda_{\infty}(\gamma)$ such that the range of $T$ and the range of $R$ is isomorphic where $R=T_{2}\circ T\circ T_{1}$. Since $T$ is tame, by using (\ref{E4}) and (\ref{E3}), it is easy to see that R is also tame. Since $R:\Lambda_{\infty}(\gamma)\to \Lambda_{\infty}(\gamma)$ is tame, the range of $R$ has a basis by Theorem 3.1 and Theorem 3.4 of \cite{vogt1}.
				Therefore the range of $T$ has a basis. 	If in addition $\Lambda_\infty(\beta)$ is nuclear, then the range of $T$ has an absolute basis since each basis in a nuclear Fr\'echet space is absolute.
			\end{proof}
			
			\begin{cor}
				Let $T\in L(\Lambda_\infty(\alpha),\Lambda_\infty(\beta))$
				be a tame operator. If $\Lambda_\infty(\beta)$ is nuclear, and the range of $T$ is closed, Then the range of $T$ is isomorphic to a closed subspace of $s$, where s is the space of rapidly decreasing sequences $\Lambda_{\infty}((\ln n)_{n\in \nc})$. 
			\end{cor}
			
			\begin{proof}
				By Proposition 31.5 in \cite{vogt4}, the range of $T$ is isomorphic to a closed subspace of $s$ since nuclearity is preserved in subspaces, 
				$\Lambda_\infty(\beta)$ has property DN, and this property is inherited by all closed subspaces.  
			\end{proof}
			
			%%%%%%%%%%%%%%%%%%%%%%%%%%%%%%%%%%%%%%%%%%%%%%%%%%%%%%%%%%%%%%%%%%%%%%%%%%%%%%%%%%%%%%%%%%%%%%%%%%%%%%%%%%%%%%%%%%%%%%%%%%%%%%%%%%%%%%%%%%%%%%%%

		\end{document}